\documentclass[a4paper]{article}
\usepackage{titlesec}
\usepackage{amsfonts,amssymb,amsmath,indentfirst,amsthm}
\usepackage{citesort}
\usepackage[sort&compress,numbers]{natbib}  
\usepackage{subfigure} 
\usepackage[footnotesize]{caption}
\usepackage[dvips]{graphics}
\usepackage{graphicx} 
\usepackage{epsfig}
\usepackage{colortbl,dcolumn}     

\begin{document}

\newtheorem{thm}{{\indent}Theorem}[section]
\newtheorem{cor}[thm]{{\indent}Corollary}
\newtheorem{lem}[thm]{{\indent}Lemma}
\newtheorem{prop}[thm]{{\indent}Proposition}
\newtheorem{result}[thm]{{\indent}Result}
\theoremstyle{definition}
\newtheorem{defn}[thm]{{\indent}Definition}
\newtheorem{rem}[thm]{{\indent}Remark}
\newtheorem{ex}[thm]{{\indent}Example}
\numberwithin{equation}{section}
\numberwithin{figure}{section}
\renewcommand{\figurename}{Fig.}
\newenvironment{keywords}{\par\textbf{keywords:}\mbox{  }}{ }
\newenvironment{Ack}{\par\textbf{Acknowledgements:}\mbox{  }}{ }
\newenvironment{coj}{\par\textbf{Conjecture:}\mbox{  }}{ }

\allowdisplaybreaks[4]

\title{Pattern formation for a volume-filling chemotaxis model with logistic growth
\thanks{The work of Yazhou Han was supported by the National Natural Science Foundation of China (Grant No. 11201443). The work of Manjun Ma Was supported by the National Natural Science Foundation of China (Grant No. 11271342), the Provincial Natural Science Foundation of Zhejiang (Grant No. LY15A010017) and  the Science Foundation of Zhejiang Sci-Tech
University under Grant No. 15062173-Y. The work of Jicheng Tao was supported by the Provincial Natural Science Foundation of Zhejiang (Grant No. LY16A010009).}}
\author{Yazhou Han$^1$\quad Zhongfang Li$^1$\quad Jicheng Tao$^1$\quad Manjun Ma$^2$
\thanks{Corresponding author, Tel: +86 571 86843224, Email: mjunm9@zstu.edu.cn }\\
\footnotesize 1. Department of Mathematics, College of
Science,\\[-0.15cm]
\footnotesize China Jiliang University, Hangzhou 310018, China\\[-0.15cm]
\footnotesize 2.Department of Mathematics, School of Science,\\[-0.15cm]
\footnotesize Zhejiang Sci-Tech University, Hangzhou, 310018, China}
\maketitle

\begin{abstract}
This paper is devoted to investigate the pattern formation of a volume-filling chemotaxis model with logistic cell growth. We first apply the local stability analysis to establish sufficient conditions of destabilization for uniform steady-state solution. Then, weakly nonlinear analysis with multi-scales is used to deal with the emerging process of patterns near the bifurcation point. For the single unstable mode case, we derive the Stuart-Landau equations describing the evolution of the amplitude,  and thus the asymptotic expressions of patterns are obtained in  both supercritical case and subcritical case. While for the case of multiple unstable modes, we also derive coupled amplitude equations to study the competitive behavior between two unstable modes through the phase plane analysis. In particular, we find that the initial data play a dominant role in the competition. All the theoretical and numerical results are in excellently qualitative agreement and better quantitative agreement than that in \cite{MOW2012}. Moreover, in the subcritical case, we confirm the existence of stationary patterns with larger amplitudes when the bifurcation parameter is less than the first bifurcation point, which gives an positive answer to the open problem proposed in \cite{MHTT2014}.


\end{abstract}

\normalsize

\section{Introduction}\label{Section_1}

In the fifties of last century, Turing \cite{T1952} had proposed a pioneering work to explain the phenomenon of pattern formation.  In \cite{T1952}, based on the reaction-diffusion system, Turing studied the self-organization process which can generate some kind of ordered structures (such as zebra stripes) in the biological world. This process is usually a passive diffusion process. Meanwhile, in the biological world, there is another class of significant diffusion --- active diffusion process, namely the tendency to some kind of chemical substances. The tendency is termed as \textit{chemotaxis}, or \textit{chemosensitive movement}. In fact, \textit{chemotaxis} plays a great role in the lives of living organisms such as locating food, searching for mates and notifying the dangers that their companions are facing.

This paper aims to study a class of chemotaxis model with volume-filling effect and logistic cell growth, which was introduced firstly by Painter and Hillen (see \cite{HP2001,PH2002}), as follows
\begin{equation}\label{eq 01}
    \begin{cases}
  \frac{\partial u}{\partial t} =\nabla(d_1\nabla u-\chi u(1-u)\nabla v)+\mu u(1-\frac{u}{u_{c}}),\\
  \frac{\partial v}{\partial t} =d_{2}\Delta v+\alpha u-\beta v,
  \end{cases}
\end{equation}
where $(x,t)\in\Omega\times[0,+\infty)$, $\Omega$ is a bounded domain in $\mathbb{R}^{N}$ with smooth boundary $\partial\Omega$; $u(x,t)$ is the cell density and $v(x,t)$ denotes the chemical concentration; $d_{1}>0$ and $d_{2}>0$ denote the cell and chemical diffusion coefficients, respectively; $\chi u(1-u)\nabla v$ represents the chemotactic flux under a volume constraint 1 (called crowding capacity), $\chi>0$ is the chemotactic coefficients measuring the strength of the chemotactic response; $\mu u(1-\frac{u}{u_{c}})$ is the cell kinetics term describing the logistic growth of cells with the growth rate $\mu>0$ and carrying capacity $u_{c}$ fulfilling $0<u_{c}<1$; $\alpha u-\beta v$ with $\alpha,\beta>0 $ is the dynamic term of chemical substances, $\alpha u$ implies that the chemical is secreted by cells themselves, $\beta v$ is the degradation of chemicals. For completeness, we shall consider the system \eqref{eq 01} subject to initial data
\begin{equation}\label{condition initial}
    u(x,0)=u_{0}\geq0,\quad v(x,0)=v_{0}\geq0,\ x\in\Omega,
\end{equation}
and Neumann boundary condition
\begin{equation}\label{condition Neumann boundary}
    \frac{\partial u}{\partial \nu}=\frac{\partial v}{\partial \nu}=0,\ t>0,x\in\partial\Omega,
\end{equation}
where $\nu$ denote the outward unit normal vector on $\partial\Omega$.

The model \eqref{eq 01}-\eqref{condition Neumann boundary} has been studied from different views by many scholars. We just outline some of them in the following paragraph. More details can be seen in \cite{HP2001, PH2002, MOW2012, MYT2010, OY2009, WX2013, MHTT2014,HLZM2016,HP2009} 
and references therein.

In 2009, Hillen and Painter\cite{HP2009} summarized the derivation of chemotaxis model and the model variations. Then, they outlined mathematical approaches for determining global existence and showed how the instability conditions depend on the parameters. Particularly, for the case $\mu= 0$, Wang and Xu\cite{WX2013} covered six specific cases of the chemotactic flux, and obtained the existence of patterns by globally bifurcation analysis. Moreover, if let $\frac\chi{d_1}\rightarrow+\infty$, Wang and Xu found that the solution of \eqref{eq 01} tended to spiky or transition layer.


In 2009, for small chemotatic parameter $\chi$, Ou and Yuan\cite{OY2009} proved the existence of travelling wavefronts connecting $(\overline{u},\overline{v}) =(u_c,\frac\alpha\beta u_c)$ to $(0,0)$, where $(u_c,\frac\alpha\beta u_c)$ and $(0,0)$ are the only two uniform steady states to \eqref{eq 01}. Ma, Yang and Tang \cite{MYT2010} generalized the result of \cite{OY2009} and established the existence of travelling wave solution for the general reaction terms. Recently, for the larger chemotatic parameter $\chi$, we studied the invasion process of the pattern in the form of a travelling wave and got the asymptotic expression of temporal-spatial pattern (see \cite{HLZM2016}).


In 2012, Ma, Ou and Wang \cite{MOW2012} studied the steady-state solutions of the volume-filling chemotaxis model \eqref{eq 01}-\eqref{condition Neumann boundary}. 
They proved that $(\bar{u},\bar{v})$ is globally asymptotically stable when $\chi=0$,  but if
\begin{equation}\label{eq 43}
    \chi>\frac{\mu d_2+\beta d_1+2\sqrt{d_1d_2\mu\beta}}{\alpha u_c(1-u_c)}\triangleq \chi_c,
\end{equation}
they gave the existence of stationary patterns by the index theory. Then, by the method of asymptotic analysis, they obtained the asymptotic expressions of stationary patterns with small amplitudes near bifurcation points and prove that, if a pattern solution with small amplitude is stable, then it must lie in the first pattern branch. However, by the method of \cite{MOW2012}, they can not obtain the specific expression of the amplitude.

In the subsequent work \cite {MHTT2014}, Ma, et al. considered the case $0<\chi<\chi_c$ and obtained the following results.

\begin{prop}\label{pro non-exis Ma}
Let $\alpha,\beta,\mu,d_1,d_2$ be fixed and take $\Omega=[0,l]$. Then, the following statements are valid:

(i) If there exists some positive integer $n_0$ such that $n_0= \left(\frac{\mu\beta}{d_1d_2} \right)^{\frac 14} \frac l\pi$, then system \eqref{eq 01} does not admit nonconstant steady-state solution with small amplitudes for any $\chi\in [0,\chi_c)$.

(ii) If $n\neq \left(\frac{\mu\beta}{d_1d_2} \right)^{\frac 14} \frac l\pi$ for any positive integer $n$, then system \eqref{eq 01} does not admit nonconstant steady-state solution with small amplitudes for any $\chi\in [0,\chi_c]$.
\end{prop}

In addition, they have proposed an open problem as follows: \textbf{whether there exist non-negative stationary pattern solutions with large amplitudes for $\chi<\chi_c$}.

In this paper, using the weakly nonlinear analysis (see \cite{GLS2012, GLS2013, GLS2014, GLSS2013,H2006,HHS1994,WMZ1994} and references therein), we will derive the amplitude equation to investigate the process of pattern formation and get much more accurate expressions of pattern solutions. Moreover, we will theoretically give a positive answer to the open problem and a stationary pattern with large amplitude will be numerically presented from the full system \eqref {eq 01}-\eqref{condition Neumann boundary} with $\chi<\chi_c$, which can be found in Subsection \ref{Sec_sub_mono}.

This paper is arranged as follows. In section \ref{Section_2}, a class of sufficient conditions of destabilization for uniform steady-state solution will be obtained by the local stability analysis. Section \ref{Section_3} is devoted to study the process of pattern formation by the weakly nonlinear analysis. We first discuss the case of single unstable mode and get the cubic and quintic Stuart-Landau equation to predict the evolution of amplitude of pattern. In particularly, under the subcritical case, we compare our results with the results of \cite{MOW2012}, which show that our  results are more efficient. Moreover, for the subcritical case, the bifurcation diagram of \eqref{eq 18} shows the phenomenon of hysteresis and jumping.
Then, the second topic discussed in Section \ref{Section_3} is the case of multiple unstable modes. There, we find that there exists competitive phenomenon between multiple unstable modes and the initial data play the dominant roles in the evolution of pattern solution. In section \ref{section_5}, we summarize the results and give some further discussion. For the completeness, four appendices, which process some specific calculation, are presented at the end of this paper.

\section{Destabilization condition}\label{Section_2}

In this section, we will give a class of sufficient conditions of destabilization for $(\overline{u},\overline{v})$ by the local stability analysis.
Let $u=\overline{u}+W^{(1)}(x,t),v=\overline{v}+W^{(2)}(x,t)$. Then, we get the linearized system of \eqref {eq 01} as follows
\begin{equation}\label{eq 02}
    \dot{\mathbf{w}}=K\mathbf{w}+D^\chi\nabla^{2}\mathbf{w},
\end{equation}
where
      $$\mathbf{w}=\begin{pmatrix} W^{(1)}\\ W^{(2)}  \end{pmatrix},\quad K=\begin{pmatrix} -\mu & 0\\ \alpha & -\beta \end{pmatrix},\quad D^\chi=\begin{pmatrix} d_{1} & -\chi u_{c}(1-u_{c}) \\ 0 & d_{2} \end{pmatrix}.$$
The solution with the form
      $\mathbf{w}=\begin{pmatrix}1\\ 1\end{pmatrix}e^{i\mathbf{k}\cdot x+\lambda t}$
leads to the following dispersion relation, which characterizes the relation between eigenvalue $\lambda$ and the wavenumber $k=|\mathbf{k}|$,
\begin{equation}\label{eq 03}
    \lambda^{2}+g(k^{2})\lambda+h(k^{2})=0,
\end{equation}
where
\begin{gather}
    g(k^{2})=k^{2}tr(D^\chi)-tr(K),\label{eq 23}\\
    h(k^{2})=det(D^\chi)k^{4}+qk^{2}+det(K),\label{eq 24}\\
    q=\mu d_{2}+\beta d_{1}-\alpha \chi u_{c}(1-u_{c}).\label{eq 25}
\end{gather}
 From the local stability theory it follows that the steady state $(\bar{u},\bar{v})$ is locally stable if $\mathbf{Re}(\lambda)\leq 0$, and $(\bar{u},\bar{v})$ is unstable if $\mathbf{Re}(\lambda)>0$.

Since
    $g(k^{2})>0,$ 
then \eqref{eq 03} does not admit a pair of conjugate pure imaginary roots. Thus, Hopf bifurcation can not appear for the system \eqref{eq 01}. 

To ensure that $(\bar{u},\bar{v})$ loses its stability, there exists at least one eigenvalue $\lambda$ satisfying $\mathbf{Re}(\lambda)>0$. Notice that 
\begin{equation}\label{eq 19}
    \lambda=\frac{1}{2}(-g\pm\sqrt{g^{2}-4h}),
\end{equation}
a class of sufficient conditions of destabilization for $(\bar{u},\bar{v})$ is $h(k^2)<0$ for some $k^2$. Next we discuss the borderline case $h(k^2)=0$.

If $h(k^2)=0$, then
\begin{equation}\label{eq 40}
    \chi=\frac{(\mu+d_1k^2)(\beta+d_2k^2)}{\alpha u_c(1-u_c)k^2}\geq 
    \chi_c,
\end{equation}
and the equal sign holds if and only if
\begin{equation}\label{eq 16}
    k^2=\sqrt{\frac{\mu\beta}{d_1d_2}}\triangleq k_c^2.
\end{equation}
On the other hand, 
if
$k^2=-\frac q{2d_1d_2},$ 
then
\begin{equation}\label{eq 20}
    h_{\text{min}}=\mu\beta-\frac{q^2}{4d_1d_2} =\mu\beta-\frac{(\alpha\chi u_c(1-u_c)-\mu d_2-\beta d_1)^2}{4d_1d_2}.
\end{equation}
Thus, as presented
\begin{figure*}[htbp]
    \centering
    \includegraphics[width=13.5cm,
    ]{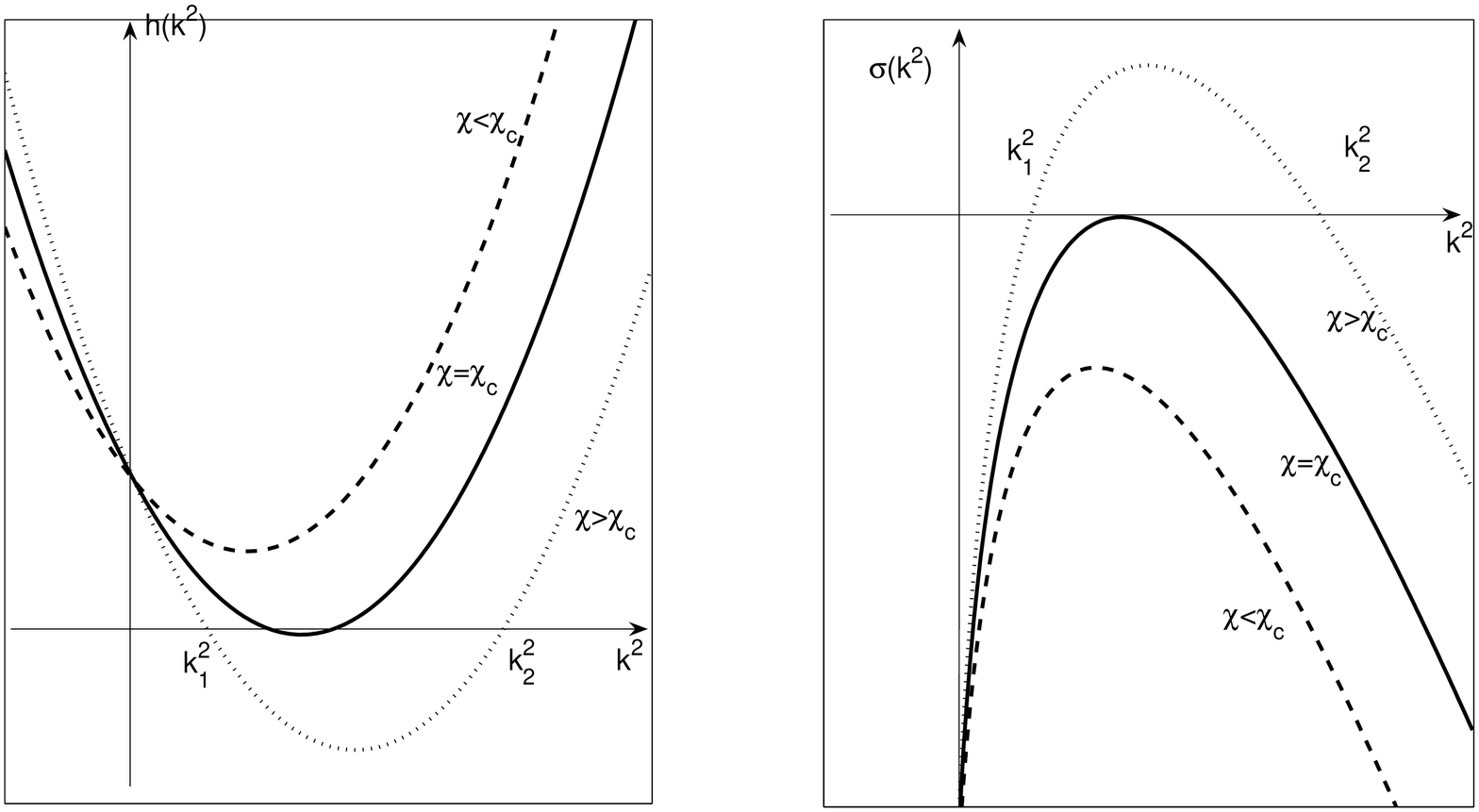}
    \caption{Left: Plot of $h(k^2)$, Right: Growth rate of the kth mode.}\label{Fig 01}
\end{figure*}
in the Fig. \ref{Fig 01},  $h_{\text{min}}=0$ provided that $\chi=\chi_c$. 
If $\chi<\chi_c$, then $h_{\text{min}}>0$ and both of the eigenvalues of \eqref{eq 03} satisfy $\mathbf{Re}(\lambda)<0$, and thus $(\bar{u},\bar{v})$ is locally asymptotically stable. If $\chi>\chi_c$, then $h_{\text{min}}<0$, and there exist $k^2$ such that \eqref{eq 03} has two real eigenvalues with different signs, which leads to the destabilization of $(\bar{u},\bar{v})$.

\begin{result}[Sufficient conditions of destabilization]\label{Result sufficient}
If $\chi>\chi_c$ and there exist modes $k^2$ such that
\begin{equation}\label{eq 47}
    k_1^2<k^2<k_2^2,
\end{equation}
then uniform steady-state solution $(\bar{u},\bar{v})$ destabilizes,
where
\begin{align}
    k_1^2=&\frac{-q-\sqrt{q^2-4d_1 d_2\mu\beta}}{2d_1d_2},\label{eq 27}\\
    k_2^2=&\frac{-q+\sqrt{q^2-4d_1 d_2\mu\beta}}{2d_1d_2}\label{eq 28}
\end{align}
are the two positive roots of $h(k^2)=0$.

In particular, when $N=1,\ \Omega=[0,l]$ , the uniform steady state $(\bar{u},\bar{v})$ destabilizes provided that $\chi>\chi_c$ and there exists at least a positive integer $n_0$ such that
   $$k_1^2<\left(\frac{n_0\pi}l\right)^2<k_2^2.$$
\end{result}

In the following, through the weakly nonlinear analysis, we shall study how the solutions evolve into stationary patterns when $\chi>\chi_c$ and $\chi$ is close to $\chi_c$. Additionally, the discuss of this paper focuses on the case of one dimension and the higher dimensional case will be presented in the subsequent work. So, we assume that $N=1$ and $\Omega=[0,l]$ in the sequel.

\section{Asymptotic expression for stationary pattern \label{Section_3}}

\subsection{Weakly nonlinear analysis}

Let
\begin{equation}\label{eq 04}
\mathbf{w}=\begin{pmatrix} W^{(1)} \\
    W^{(2)} \end{pmatrix}=\begin{pmatrix} u-\overline{u} \\ v-\overline{v} \end{pmatrix},
\end{equation}
and then the system \eqref{eq 01} can be rewritten as follows
\begin{equation}\label{eq 05}
    \mathbf{w_{t}}=\mathcal{L}^{\chi}\mathbf{w} +\mathcal{N}\mathcal{L}^{\chi}\mathbf{w},
\end{equation}
where
    $$\mathcal{L}^{\chi}=K+D^{\chi} \frac{\partial^2}{\partial x^2},$$
    $$\mathcal{N}\mathcal{L}^{\chi}\mathbf{w}= \begin{pmatrix} -\chi\frac\partial{\partial x}\bigl((1-2u_{c})W^{(1)}\frac{\partial W^{(2)}}{\partial x} -(W^{(1)})^{2}\frac{\partial W^{(2)}}{\partial x} \bigr)-\frac\mu{u_{c}} {(W^{(1)})^{2}} \\ 0 \end{pmatrix}.$$
Expand $\chi,\mathbf{w}$ and $t$ as
\begin{align}
    &\chi=\chi_{c} +\varepsilon\chi_{1} +\varepsilon^{2}\chi_{2} +\varepsilon^{3}\chi_{3} +\varepsilon^{4}\chi_{4} +\varepsilon^{5}\chi_5 +\cdots,\label{eq 29}\\
    &\mathbf{w}=\varepsilon\mathbf{w_{1}} +\varepsilon^{2}\mathbf{w_{2}} +\varepsilon^{3}\mathbf{w_{3}} +\varepsilon^{4}\mathbf{w_{4}} +\varepsilon^{5}\mathbf{w_5}+\cdots,\label{eq 30}\\
    &t=t(T_1,T_2,T_3,\cdots),\quad T_i=\varepsilon^i t, i=1,2,\cdots,\label{eq 31}
\end{align}
where $\mathbf{w}_i=\bigl(W_i^{(1)},W_i^{(2)}\bigr)^T$, and $T_i, i=1,2,\cdots$ represent different time scales.
Substituting \eqref{eq 29}, \eqref{eq 30} and \eqref{eq 31} into \eqref{eq 05} and collecting the terms at each order in $\varepsilon$, we obtain a sequence of coefficient equations.

\noindent $O(\varepsilon)$:
\begin{equation}\label{eq 06}
\mathcal{L}^{\chi_{c}}\mathbf{w_{1}}=0,
\end{equation}
$O(\varepsilon^{2})$:
\begin{equation}\label{eq 07}
\mathcal{L}^{\chi_{c}}\mathbf{w_{2}}=\mathbf{F} =\begin{pmatrix}F^{(1)}\\ F^{(2)} \end{pmatrix} =\begin{pmatrix}F^{(1)}\\ \frac{\partial W_1^{(2)}}{\partial T_1}\end{pmatrix},
\end{equation}
$O(\varepsilon^{3})$:
\begin{equation}\label{eq 08}
\mathcal{L}^{\chi_{c}}\mathbf{w_{3}}=\mathbf{G},
\end{equation}
where
\begin{equation}\label{eq 32}
    \begin{split}
    F^{(1)}=&\frac{\partial W_{1}^{(1)}}{\partial T_{1}}+\chi_1 u_c(1-u_c)\frac{\partial^2 W_{1}^{(2)}}{\partial x^2}\\ &+\chi_{c}\frac\partial{\partial x}((1-2u_{c})W_{1}^{(1)}\frac\partial{\partial x} W^{(2)}_{1}) +\frac{\mu}{u_{c}} (W^{(1)}_{1})^{2}.
    \end{split}
\end{equation}
The explicit expression of $G$ and all the detailed calculation are given in Appendix \ref{appendix coefficient-eqs}.

\subsection{Stationary pattern for single unstable mode}

Firstly, we assume that $k_c^2$ is the unique unstable mode satisfying \eqref{eq 47}. So, the solution to the equation \eqref{eq 06} with the Neumann boundary conditions is
\begin{equation}\label{eq 09}
    \mathbf{w_{1}}=A\mathbf{\rho}\cos(k_{c}x),
\end{equation}
where the amplitude function $A$ only depends on temporal variable and
\begin{equation}\label{eq 10}
    \mathbf{\rho}=\begin{pmatrix} M\\ 1 \end{pmatrix}=\begin{pmatrix} \frac{\beta+k_{c}^2d_{2}}{\alpha}\\ 1 \end{pmatrix}\in \text{Ker}(K-k_c^2D^\chi).
\end{equation}
Substituting $\mathbf{w_{1}}$ into  \eqref{eq 07} 
leads to
\begin{equation}\label{eq 14}\begin{split}
    \mathbf{F}&=\left[\frac{\partial A}{\partial T_1}\rho- A\begin{pmatrix} 0 & k_{c}^{2}\chi_{1}u_c(1-u_c)\\ 0& 0\\ \end{pmatrix} \rho\right]\cos(k_{c}x)\\
    &- \begin{pmatrix} M\chi_{c}k_{c}^{2}(1-2u_{c})-\frac{\mu M^2}{2u_c}\\ 0 \end{pmatrix} A^{2}\cos(2k_{c}x)+ \begin{pmatrix} \frac{\mu M^2}{2u_{c}}\\ 0 \end{pmatrix}A^{2}.
\end{split}\end{equation}
Suppose that
\begin{equation}\label{eq 11}
    \mathbf{w}^*=\psi\cos(k_cx)\quad \text{with}\ \psi=\begin{pmatrix} M^*\\ 1 \end{pmatrix},
\end{equation}
is a fundamental solution of
    $L^*\mathbf{w}^*=0,$
where $L^*$ is the adjoint operator of $L^{\chi_c}$, $M^*=\frac\alpha{\mu+d_1 k_c^2}$.\ By the solvability condition for \eqref{eq 07}, we have $<\mathbf{F},\mathbf{w}^*>=0$,
and then obtain
\begin{equation}\label{eq 12}
    \frac{\partial A}{\partial T_1}=\gamma A\quad\text{with } \gamma=\frac{k_{c}^{2}\chi_{1}(u_{c}-u_{c}^{2})M^*} {1+MM^{*}}.
\end{equation}
Noting that the solution $A=C\exp(\gamma T_1)$ of \eqref{eq 12} can not predict correctly the evolution of amplitude, we take $\frac{\partial A}{\partial T_{1}}=0$ and $\chi_{1}=0$ simply to satisfy the solvability condition. Particularly, $\frac{\partial A}{\partial T_{1}}=0$ means that solution is independent of time scale $T_1$, i.e. $\frac{\partial\mathbf{w}}{\partial T_1}=0$.

Then, from \eqref{eq 07} and \eqref{eq 14} it follows that \eqref{eq 07} has solution
\begin{equation}\label{eq 33}
    \mathbf{w_{2}}=A^{2}\mathbf{W_{20}} +A^{2}\mathbf{W_{22}}\cos(2k_{c}x),
\end{equation}
where
\begin{equation}\label{eq 34}
    \mathbf{W_{20}}= \begin{pmatrix} -\frac{M^2}{2u_c}\\ -\frac{\alpha M^2}{2\beta u_c} \end{pmatrix},\quad \mathbf{W_{22}}= \begin{pmatrix} \frac{\beta+4k_c^2 d_2}\alpha \\ 1\end{pmatrix}\mathbf{W}_{22}^{(2)},
\end{equation}
\begin{equation}\label{eq 35}
    \mathbf{W}_{22}^{(2)}= \frac{\alpha\Bigl(\frac{\mu M^2}{2u_{c}}-\chi_{c}k_{c}^{2}M(1-2u_{c})\Bigr)}{4\alpha k_c^2\chi_cu_c(1-u_c)-(\mu+4k_c^2 d_1)(\beta+4k_c^2 d_2)}.
\end{equation}
To explore the evolution of the amplitude $A$, we further discuss the third-order coefficient equation \eqref{eq 08}. Substituting $\mathbf{w}_1$ and $\mathbf{w}_2$ into \eqref{eq 08} and combining with the solvability condition  $<\mathbf{G},\mathbf{w}^*>=0$, we get the cubic Stuart-Landau equation of amplitude $A$ as follows
\begin{equation}\label{eq 13}
    \frac{d A}{d T_{2}}=\sigma A- LA^{3},
\end{equation}
where
\begin{align}
    &\sigma=\frac{\langle G_{11},\psi\rangle}{\langle \rho,\psi\rangle} =\frac{k_{c}^{2}\chi_{2}(u_{c}-u_{c}^{2})M^*} {1+MM^*}>0,\label{eq 37}\\
     &L=\frac{\langle G_{13},\psi\rangle}{\langle \rho,\psi\rangle}=\frac{G_{13}^{(1)}M^*} {1+MM^*},\label{eq 38}\\
     &G_{13}^{(1)}=\chi_{c}(2u_{c}-1)k_{c}^{2}\bigl(W^{(1)}_{20}+MW^{(2)}_{22} -\frac 12 W^{(1)}_{22}\bigr)\nonumber\\
     &\hspace{1.2cm}+\frac{1}{4}\chi_{c}M^2k_{c}^{2} +\frac{\mu M}{u_{c}}(2W^{(1)}_{20}+W^{(1)}_{22}).\label{eq 39}
\end{align}
The detailed computations of \eqref{eq 33} and \eqref{eq 13} can be found in Appendix \ref{appendix Stuart-Landau}.

Obviously, the dynamics of \eqref{eq 13} can be divided into two cases according to the sign of $L$, i.e., the supercritical case and the subcritical case corresponding to $L>0$ and $L<0$, respectively. By \eqref{eq 38}, the sign of $L$ completely depends on that of $G_{13}^{(1)}$. Especially, when $u_c=\frac 12$, we have
\begin{equation}\label{eq 26}\begin{split}
    G_{13}^{(1)}=&\frac 14\chi_cM^2k_c^2-4\mu M^3\\
    &+\frac{2\mu^2 M^3(\beta+4k_c^2 d_2)} {4k_c^2\alpha \chi_c u_c(1-u_c) -(\mu+4k_c^2d_1)(\beta+4k_c^2d_2)}\\
    =&\sqrt{\mu}\left\{\frac 14\chi_cM^2\sqrt{\frac\beta{d_1d_2}} -4\sqrt{\mu} M^3\right.\\
    &\left.\hspace{1.0cm}+\frac{2\mu^{\frac 32}M^3(\beta+4k_c^2 d_2)} {k_c^2\alpha \chi_c -(\mu+4k_c^2d_1)(\beta+4k_c^2d_2)}\right\}.
\end{split}\end{equation}
Thus, we can obtain the following conclusion.
\begin{result}\label{Result sign L with uc}
 Let $u_c=\frac 12$ and other parameters $d_1, \ d_2,\ \alpha$ and $\beta$ be fixed. Then, if the cell growth rate $\mu$ is small enough, then \eqref{eq 13} is supercritical; while $\mu$ is large enough, then \eqref{eq 13} is subcritical.
\end{result}
In the following, we shall derive the amplitude equation in both supercritical case and subcritical bifurcation case.

\subsubsection{The supercritical case}\label{Sec_sup_mono}

\begin{figure*}[hbpt]
\centering
\includegraphics[width=13.5cm
]{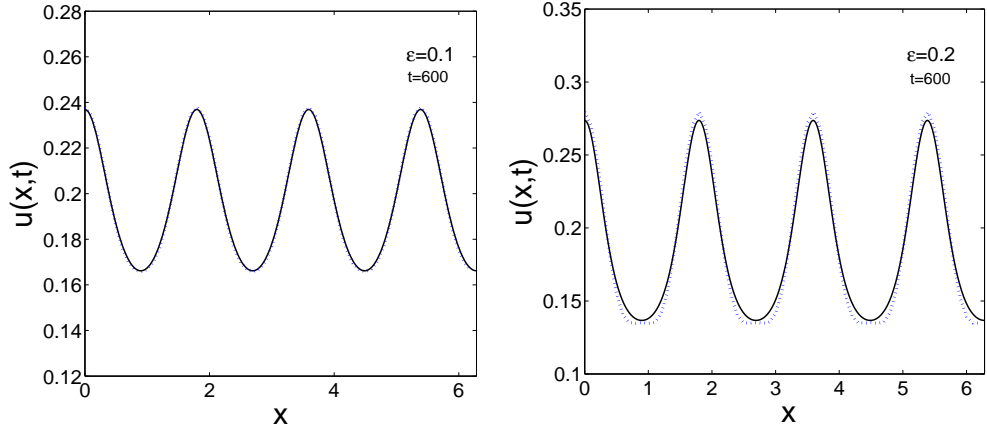}
\caption{Comparison between the weakly nonlinear solution (dashed line) and the numerical solution of \eqref{eq 01}-\eqref{condition Neumann boundary}(solid line). The initial data is set as $1\%$ random small perturbation of the $(u_c,\frac\alpha\beta u_c)$ and the parameters are $\alpha=36,\beta=34,\mu=0.5,u_{c}=0.2,d_{1}=0.2,d_{2}=0.6, l=2\pi$. With these parameters one has $\chi_{c}= 1.7286,k_c= 3.45$ and the first admissible unstable mode is $\overline{k_{c}}=3.5$.}\label{Fig 05}
\end{figure*}

Because of $\sigma>0$ and $L>0$, it is easy to know that \eqref{eq 13} has a globally asymptotic stable solution 
    $$A_{\infty}\triangleq\sqrt{\cfrac{\sigma}{L}},$$
which represents the limit value of the amplitude $A$. Substituting $A_\infty$ into \eqref{eq 30}, \eqref{eq 09} and \eqref{eq 33}, we have the second-order asymptotic expression of the stationary pattern as follows
\begin{equation}\label{eq 17}
    \mathbf{w}=\varepsilon\rho\sqrt{\frac{\sigma}{L}}\cos(k_{c}x) +\varepsilon^{2}\frac{\sigma}{L}(\mathbf{w_{20}} +\mathbf{w_{22}}\cos(2k_{c}x))+O(\varepsilon^{3}).
\end{equation}
In Fig.\ref{Fig 05},
a numerical example is presented to show the comparison between the numerical solution of \eqref{eq 01}-\eqref{condition Neumann boundary} and the weakly nonlinear asymptotic solution \eqref{eq 17}. Since $k_c\neq \frac{n\pi}l$ for any positive integer $n$, in \eqref{eq 17} we replace $k_c$ by $\overline{k_{c}}=3.5$, i.e., the first unstable mode when $\chi$ passes the critical value $\chi_c$. It is observed from Fig.\ref{Fig 05} that both the wave number and the amplitude are in perfectly agreement between the numerical solution of \eqref{eq 01}-\eqref{condition Neumann boundary} and the approximation \eqref{eq 17}  for the cases $\varepsilon=0.1$ and $\varepsilon=0.2$.

\begin{rem}\label{rem comp_2order}
If there is no $n\in N_+$ such that $k_c=\frac{n\pi}l$ for $\chi=\chi_c$, then we replace $k_c$ by the first admissible mode $\overline{k_c}=\frac{n_0\pi}l$ in \eqref{eq 09} and \eqref{eq 33}, where $n_0$ satisfies $\chi_0(n_0)= \chi_{\text{min}}$ with
    $$\chi_{\text{min}}= \min_{n}\left\{\chi_0(n)=\frac{\bigl(d_1(\frac{n\pi}l)^2+\mu\bigr) \bigl(\beta+d_2(\frac{n\pi}l)^2\bigr)} {\alpha u_c(1-u_c)(\frac{n\pi}l)^2},\ n=1,2,\cdots\right\}.$$
Here $\chi_{\text{min}}$ is the same as that in page 753 of \cite{MOW2012}.
Accordingly, the critical value $\chi_c$ is replaced by the first bifurcation value $\chi_{\text{min}}$. On the other hand, we note that the second-order weakly nonlinear asymptotic solution \eqref{eq 17} is consistent with that obtained in \cite{MOW2012}. But in \cite{MOW2012} the amplitude was not explicitly expressed by the cubic Stuart-Landau equation \eqref{eq 13}. So we improve the related results to the reference \cite{MOW2012}.
\end{rem}
\begin{rem}\label{rem mu=0}
If $\mu=0$, then
    $$k_c=\frac{n_0\pi}l=0$$
with $n_0=0$. Obviously there is no pattern formation.
Thus, when $\mu=0$, the bifurcation parameter $\chi$ should start with $\chi_{\text{min}}$.
\end{rem}

\subsubsection{The subcritical case}\label{Sec_sub_mono}

If $\sigma>0$ and $L<0$, it is easy to know that cubic Stuart-Landau equation \eqref{eq 13} is not able to capture the amplitude of the pattern. Therefore, we need to push the weakly nonlinear expansion to higher order.

Performing the weakly nonlinear expansion up to $O(\varepsilon^5)$, we can get the quintic Stuart-Landau equation for the amplitude as follows
\begin{equation}\label{eq 18}
    \frac{dA}{dT}=\bar{\sigma}A-\bar{L}A^3+\bar{Q}A^5,
\end{equation}
where
\begin{equation}\label{eq 45}
    \bar{\sigma}=\sigma+\varepsilon^{2}\widetilde{\sigma}, \quad \bar{L}=L+\varepsilon^{2}\widetilde{L}, \quad \bar{Q}=\varepsilon^{2}\widetilde{Q},
\end{equation}
$\sigma$ and $L$ are given in\eqref{eq 37} and \eqref{eq 38}, respectively; and $\tilde{\sigma}, \tilde{L}$ and $\tilde{Q}$ are given in \eqref{eq A27}. The detailed calculations are given in Appendix \ref{appendix Stuart-Landau 5}.

In the subcritical case, namely $\bar{\sigma}>0$ and $\bar{L}<0$, and when $\bar{Q}<0$, it is easy to check that
\eqref{eq 18} has a globally asymptotic stable solution
    $$\bar{A}_\infty=\sqrt{\frac{\bar{L} -\sqrt{\bar{L}^2 -4\bar{\sigma}\bar{Q}}} {2\bar{Q}}},$$
which is the limit value of amplitude $A$. Substituting $\bar{A}_\infty$ into $\mathbf{w}_1, \mathbf{w}_2, \mathbf{w}_3, \mathbf{w}_4$, we obtain the following fourth-order weakly nonlinear asymptotic expression of the stationary pattern
\begin{equation}\label{eq 46}
    \begin{split}
    \mathbf{w}=&\varepsilon\mathbf{w}_1 +\varepsilon^2\mathbf{w}_2 +\varepsilon^3\mathbf{w}_3 +\varepsilon^4\mathbf{w}_4 +O(\varepsilon^5)\\
    =&\varepsilon\bar{A}_\infty\rho\cos(k_cx) +\varepsilon^2\bar{A}_\infty^2\left(\mathbf{W}_{20} +\mathbf{W}_{22}\cos(2k_cx)\right)\\
    &+\varepsilon^3\left[\bigl(\bar{A}_\infty \mathbf{W}_{31}+\bar{A}_\infty^3 \mathbf{W}_{32}\bigr)\cos(k_cx) +\bar{A}_\infty^3 \mathbf{W}_{33} \cos(3k_cx)\right]\\
    &+\varepsilon^4\left[\bar{A}_\infty^2 \mathbf{W}_{40}+\bar{A}_\infty^4 \mathbf{W}_{41}+ \bigl(\bar{A}_\infty^2 \mathbf{W}_{42} +\bar{A}_\infty^4 \mathbf{W}_{43}\bigr) \cos(2k_cx)\right.\\
    &\left.+\bar{A}_\infty^4 \mathbf{W}_{44} \cos(4k_cx)\right]+O(\varepsilon^5),
    \end{split}
\end{equation}
where $\rho$  is given in \eqref{eq 10}, $\mathbf{W}_{20}$ and $\mathbf{W}_{22}$  are given in \eqref{eq 34}, $\mathbf{W}_{31}$, $\mathbf{W}_{32}$ and $\mathbf{W}_{33}$  are given in \eqref{eq A23}, and $\mathbf{W}_{40}$, $\mathbf{W}_{41}$, $\mathbf{W}_{42}$, $\mathbf{W}_{43}$ and $\mathbf{W}_{43}$  are given in \eqref{eq A25}.

\begin{figure*}[hbpt]
\centering
\includegraphics[width=13.5cm
]{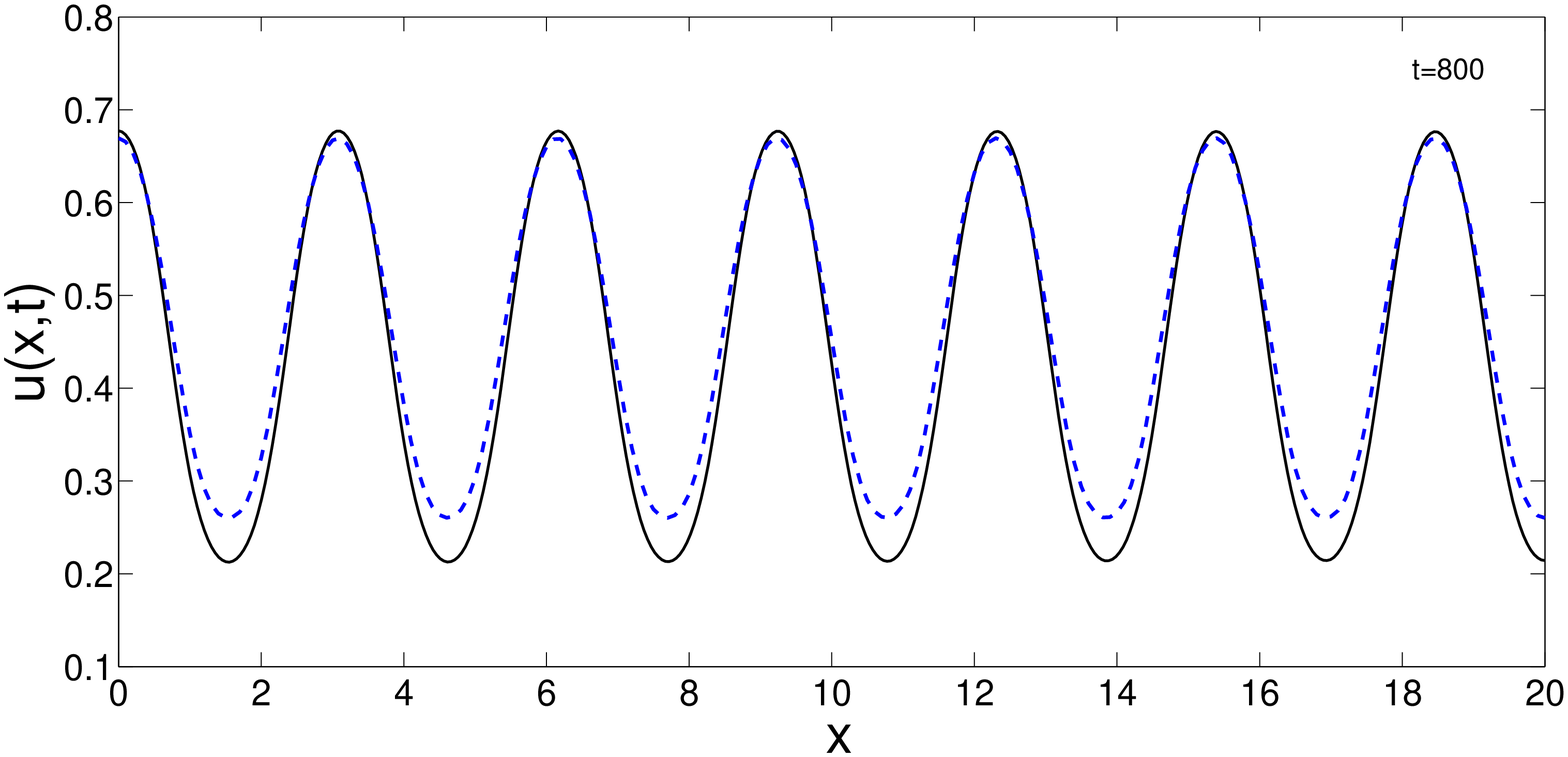}
\caption{Comparison between the weakly nonlinear solution (dashed line) and the numerical solution of \eqref{eq 01}-\eqref{condition Neumann boundary}(solid line) with $\varepsilon=0.1$, where the initial data is same as in Fig. \ref{Fig 05} and the parameters are $\alpha=10,\ \beta=10,\ \mu=0.5,\ u_{c}=0.5,\ d_{1}=0.3,\ d_{2}=1,\ l=20$. In this case, $\chi_{c}= 2.3798,\ k_c= 2.0205$.}\label{Fig 08}
\end{figure*}
Fig.\ref{Fig 08} shares the same parameters as that in Fig. 4 of \cite{MOW2012}. We have $\bar{\sigma}=1.5351,\ \bar{L}=-0.7588$ and $\bar{Q}=-0.6415$. And then, the equation \eqref{eq 18} has a stable equilibrium $\bar{A}_\infty=1.4992.$  we see that in Fig.\ref{Fig 08} the analytical approximation and the numerical simulation are in qualitative and  quantitative agreement better than Fig.4 of \cite{MOW2012}.

\begin{figure*}[hbpt]
\centering
\includegraphics[width=13.5cm
]{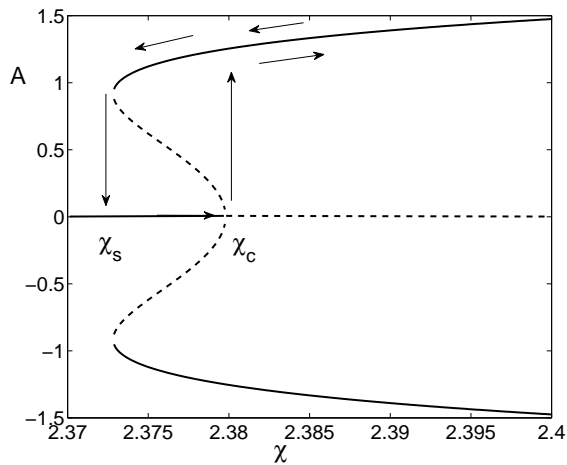}
\caption{The bifurcation diagram of \eqref{eq 18}. The parameter values are the same as that in Figure \ref{Fig 08}}\label{Fig 14}
\end{figure*}
With the parameters in Fig.\ref{Fig 08}, we plot the bifurcation diagram Fig.\ref{Fig 14}
of \eqref{eq 18}. In Fig. \ref{Fig 14}, there exist two turning points $\chi_c$ and $\chi_s$, where $\chi_s$  is the zero point of the discriminant of the equilibrium equation of \eqref{eq 18}, i.e. $\chi_s$ is the root of
    $\bar{L}^2-4\bar{\sigma}\bar{Q}=0.$
In figure \ref{Fig 14}, we numerically get
      $$\chi_c= 2.3798,\quad \chi_s= 2.3728.$$
Specifically, Fig.\ref{Fig 14} implies the following information. the origin is locally stable for $\chi<\chi_c$, and when $\chi=\chi_c$, two backward unstable branches are generated from the origin. Decreasing $\chi$ until $\chi=\chi_s<\chi_c$, these unstable branches turn around and become stable.  Thus,  two stable branches coexists in the range $\chi_s<\chi<\chi_c$.

The coexisting phenomenon indicates that solutions will be hysteresis in the range $\chi_s<\chi<\chi_c$ and jump rapidly at the two values $\chi=\chi_s$ and $\chi=\chi_c$. Specifically, for any given small perturbation around $(u_c, \frac{\alpha}{\beta}u_c)$, the solution approaches asymptotically to $(u_c, \frac{\alpha}{\beta}u_c)$ if $\chi<\chi_c$; while for $\chi>\chi_c$, the solution will jump rapidly to the stable equilibrium with large amplitude. Similarly, for the stable branch with larger amplitude, there exist the hysteresis and the jumping phenomenon at the value $\chi=\chi_s$.

\begin{figure*}[hbpt]
\centering
\includegraphics[width=13.5cm,
]{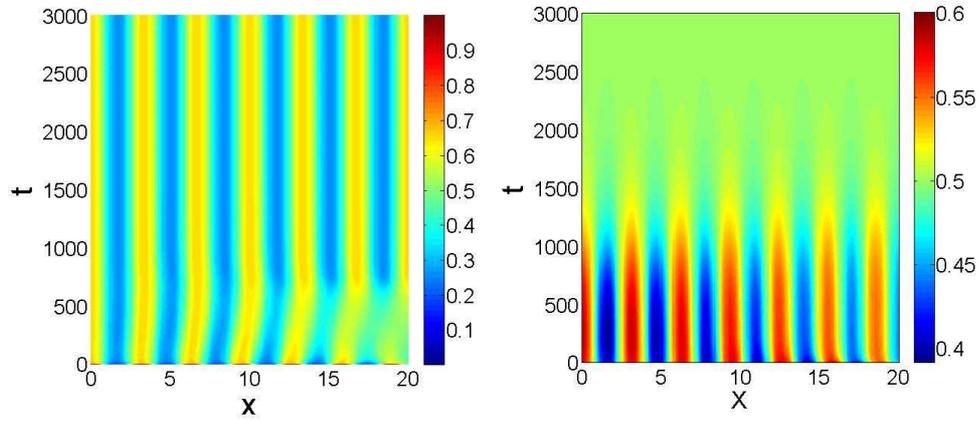}
\caption{two stable steady states $\mathbf{w}$ expressed by \eqref{eq 46} and $(u_c, \frac{\alpha}{\beta}u_c)$ coexist.   The parameters are the same as that in Fig. \ref{Fig 08} and $\chi=2.3760\in(\chi_s,\chi_c)$. Left: Pattern $\mathbf{w}$ reached by giving the initial value $u_0=u_c+0.5\cos(2x)$. Right: the uniform steady state $(u_c, \frac{\alpha}{\beta}u_c)$ reached by giving the initial value $u_0=u_c+0.1\cos(2x)$.
 }
\label{Fig 15}
\end{figure*}
Moreover, the coexisting phenomenon implies that there are stationary patterns with larger amplitudes when the bifurcation parameter $\chi\in(\chi_s,\chi_c)$. This gives an affirmative answer to the open problem  proposed in \cite{MHTT2014}. A numerical example is presented in Fig. \ref{Fig 15}.

\subsection{Stationary pattern for double unstable modes}\label{Sec_sup_multi}

\begin{figure*}[hbpt]
\centering
\includegraphics[width=13.5cm
]{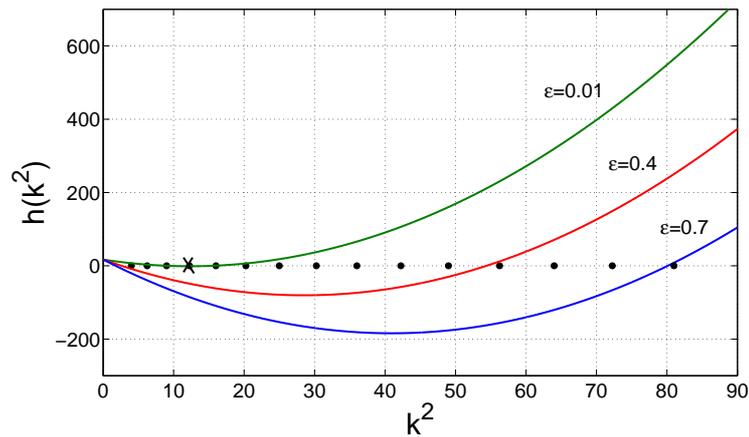}
\caption{The unstable modes are the integers and semi-integers (marked with dots) for which $h(k^2)<0$. The critical mode $\bar{k_c}$ is marked with a cross. The parameters are same as Fig. \ref{Fig 05}.} \label{Fig 07}
\end{figure*}

As presented in Fig. \ref{Fig 07}, when $\varepsilon$ is small, there exists only one mode such that $h(k^2)<0$. Then, $\overline{k_c}$ is just the most unstable mode. For this case, the asymptotic expressions of stationary pattern are discussed in the above Section. While, for larger $\varepsilon$, the parameter $\chi$ has a larger deviation from $\chi_c$, there will be more than one modes satisfy $h(k^2)<0$. We shall study how the unstable modes interact and how to determine the shape of stationary pattern.

\begin{figure*}[hbpt]
\centering
\includegraphics[width=13.5cm
]{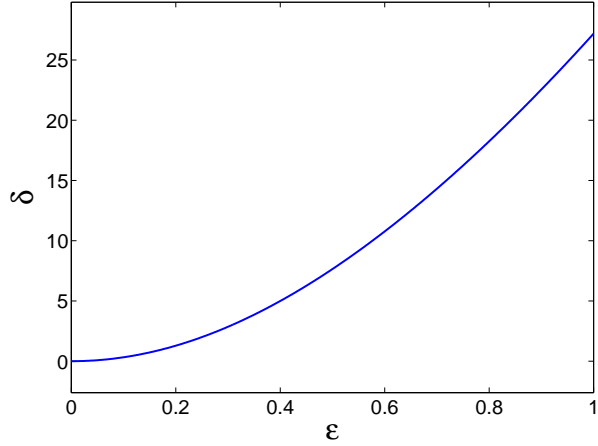}
\caption{The curve of $\frac{d\lambda^+}{d(k^2)}=0$, where the parameters are chosen as in Fig. \ref{Fig 05}.} \label{Fig 03}
\end{figure*}

We first analyze how the most unstable mode change. Denote the  positive eigenvalue of \eqref{eq 19} as $\lambda^+$, and then
    $$\lambda^+(k^2)=\frac{1}{2}(-g+\sqrt{g^{2}-4h}).$$
Combining the method in \cite{S1987}, we know that the solution $k_m^2$ of
\begin{equation}\label{eq 21}
    \frac{d\lambda^+}{d(k^2)}=0
\end{equation}
maximizes the growth rate $\lambda^+$. Thus $k_m$ is the wave number of the most unstable mode. Substituting
    $\chi=\chi_c(1+\varepsilon^2)\quad \text{and}\quad k_m^2=k_c^2+\delta$
into \eqref{eq 21}, we have
\begin{equation}\label{eq 41}
    \delta=\frac{-B-\sqrt{B^2-4AC}}{2A}-k_c^2,
\end{equation}
where
\begin{gather*}
    A=-(d_1-d_2)^2d_1d_2<0,\\
    B=4d_1d_2q-2d_1d_2(d_1+d_2)(\mu+\beta)<0,\\
    C=q^2-(d_1+d_2)(\mu+\beta)q+(d_1+d_2)^2\mu\beta>0,\\
    q=-\varepsilon^2(\mu d_2+\beta d_1)-2(1+\varepsilon^2)\sqrt{\mu\beta d_1d_2}<0.
\end{gather*}
Noting that $\delta=0$ if $\varepsilon=0$ and $\frac{d\delta}{d(\varepsilon^2)}>0$ if $\varepsilon>0$, so we know that, with the increasing of $\chi$, the most unstable mode $k_m^2$ will move away from $k_c^2$. Numerical examples are given in Fig. \ref{Fig 03}.

Then, we will investigate the competitive law between two unstable modes  $k_{1}$  and  $k_{2}$ by deriving their amplitude equations. %
Set the solution of \eqref{eq 06} as
\begin{equation}\label{eq 15}
    \mathbf{w_{1}}=\mathop{\sum}^{2}_{l=1} A_{l}\rho_{l}\cos(k_{l}x),
\end{equation}
where  $A_l$  depending only on temporal variable is the amplitude of mode $k_l$ and
    $$\rho_{l}= \begin{pmatrix}M_{l} \\ 1\end{pmatrix}\quad \text{with } M_{l}=\frac{\beta+k_{l}^{2}d_{2}}{\alpha},\  l=1,2.$$
Applying the Fredholm alternative to \eqref{eq 07} and  \eqref{eq 08} and repeating the process establishing the amplitude equation \eqref{eq 13}, we obtain the following ODE model
\begin{equation}\label{eq 22}
    \begin{cases}
        \frac{d A_{1}}{d T} =\sigma_{1}A_{1} -L_{1}A_{1}^{3} -\Omega_{1}A_{1}A_{2}^{2},\\
        \frac{d A_{2}}{d T} =\sigma_{2}A_{2} -L_{2}A_{2}^{3} -\Omega_{2}A_{2}A_{1}^{2}.
    \end{cases}
\end{equation}
The detailed calculations are given in Appendix \ref{appendix competing-eqs}. Therefore, the first-order asymptotic expression of the stationary pattern is as follows
\begin{equation}\label{eq 48}
    \mathbf{w}=\varepsilon (\rho_1 A_{1\infty}\cos(k_1 x)+\rho_2 A_{2\infty}\cos(k_2 x_)+O(\varepsilon^2),
\end{equation}
where $(A_{1\infty},A_{2\infty})$ is some stable stationary state of the system \eqref{eq 22}.

Since a complete analysis of the stationary points of \eqref{eq 22} is too involved, therefore we just present a numerical study of a typical case. Now, the values of parameters are chosen the same as in Fig. \ref{Fig 05} except $\varepsilon=0.4$.  Then, the solution of \eqref{eq 21} is about $k_m^2=(4.1105)^2$. So, we choose two unstable modes $k_1=4$ and $k_2=3.5$. By Appendix \ref{appendix competing-eqs}, 
we have
\begin{equation*}\begin{matrix}
    \sigma_1 = 3.3680,   & L_1 = 41.2467, & \Omega_1 = 75.1183,\\
    \sigma_2 = 2.7532,   & L_2 = 28.1224, & \Omega_2 = 62.0933
\end{matrix}\end{equation*}
and the four non-negative equilibria
\begin{equation*}
    (0,0),\;\;(0,0.3129),\;\;(0.2858,0),\;\;(0.1995,0.1475).
\end{equation*}
A further calculation leads to the corresponding Jacobian matrices as follows
\begin{gather*}
    J_{(0,0)}=\begin{pmatrix} 3.3680 & 0\\ 0 & 2.7532\end{pmatrix},\ J_{(0,0.3129)}=\begin{pmatrix} -3.9862 & 0\\ 0 &-5.5065\end{pmatrix},\\
    J_{(0.2858,0)}=\begin{pmatrix} -6.7359& 0\\ 0 & -2.3170\end{pmatrix},\ J_{(0.1995,0.1475)}=\begin{pmatrix} -3.1898& -4.4202\\ -3.6538 & -1.5529\end{pmatrix}.
\end{gather*}
Therefore,
\begin{figure*}[hbpt]
\centering
\includegraphics[width=10cm
]{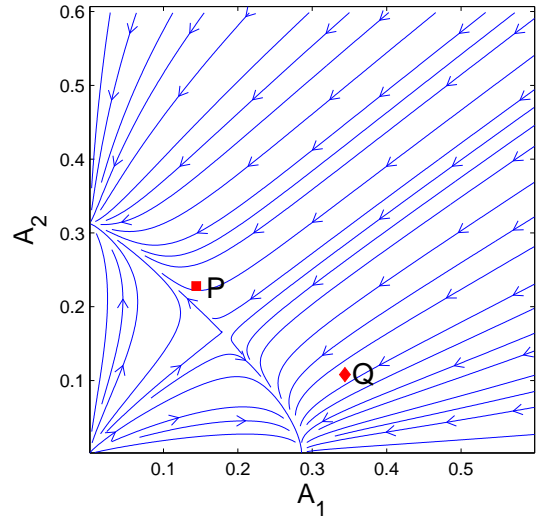}
\caption{The phase portrait of the competitive model \eqref{eq 22} for modes $k_1=4$ and $k_2=3.5$. The parameters are the same as that in Fig. \ref{Fig 05} except $\varepsilon=0.4$.}\label{Fig 11}
\end{figure*}
$(0,0)$ is an unstable node, $(0.1995,0.1475)$ is a saddle point, and $(0,0.3129)$ and $(0.2858,0)$ are stable nodes. These can also be observed in the phase diagram in Fig.\ref{Fig 11}. So, the model \eqref{eq 22} with these values is competitive and the system can evolve toward one of the two semi-trivial stable states $(A_1,0)=(0.2858,0)$ and $(0,A_2)=(0,0.3129)$. If we take the initial data $(0.144,0.228)$, which corresponds to the point $P$ in Fig \ref{Fig 11} and lies in the basin of attraction of the equilibrium $(0,A_2)$, then the asymptotic expression of stationary pattern is
\begin{equation}\label{eq 49}
    u=u_{c} +\varepsilon{A_{2}}M_2\cos(k_{2}x) +O(\varepsilon^2).
\end{equation}
The detailed comparison between numerical solution of \eqref{eq 01}-\eqref{condition Neumann boundary} and expected solution \eqref{eq 49} is presented in Fig \ref{Fig 09}.
\begin{figure*}[hbpt]
\centering
\includegraphics[width=13.5cm
]{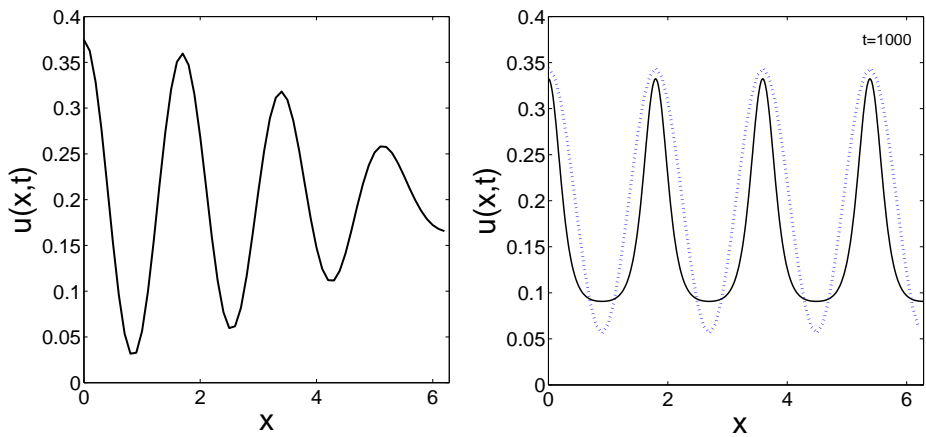}
\caption{Left:Initial condition $u=u_{c}+ \varepsilon\bigl(\bar{A_{1}}M_1\cos(k_{1}x)+ \bar{A_{2}}M_2\cos(k_{2}x)\bigr)$ with $k_{1}=4$ and $k_{2}=3.5$, where $\bar{A_{1}}=0.144,\bar{A_{2}}=0.228$ is in the basin of attraction of the equilibrium $(0,A_{2})$. Right: The comparison between the numerical solution of system \eqref{eq 01}-\eqref{condition Neumann boundary} (solid line) and the expected solution \eqref{eq 49}.}\label{Fig 09}
\end{figure*}
While for the initial data $Q(0.344,0.108)$ belonging the basin of attraction of the equilibrium $(A_1,0)$ (see Fig \ref{Fig 11}), the asymptotic expression of the stationary pattern is changed into
\begin{equation}\label{eq 50}
    u=u_{c} +\varepsilon{A_{1}}M_1\cos(k_{1}x) +O(\varepsilon^2).
\end{equation}
The comparison between \eqref{eq 50} and the simulation solution of \eqref{eq 01}-\eqref{condition Neumann boundary} is given in Fig \ref{Fig 10}.
\begin{figure*}[hbpt]
\centering
\includegraphics[width=13.5cm
]{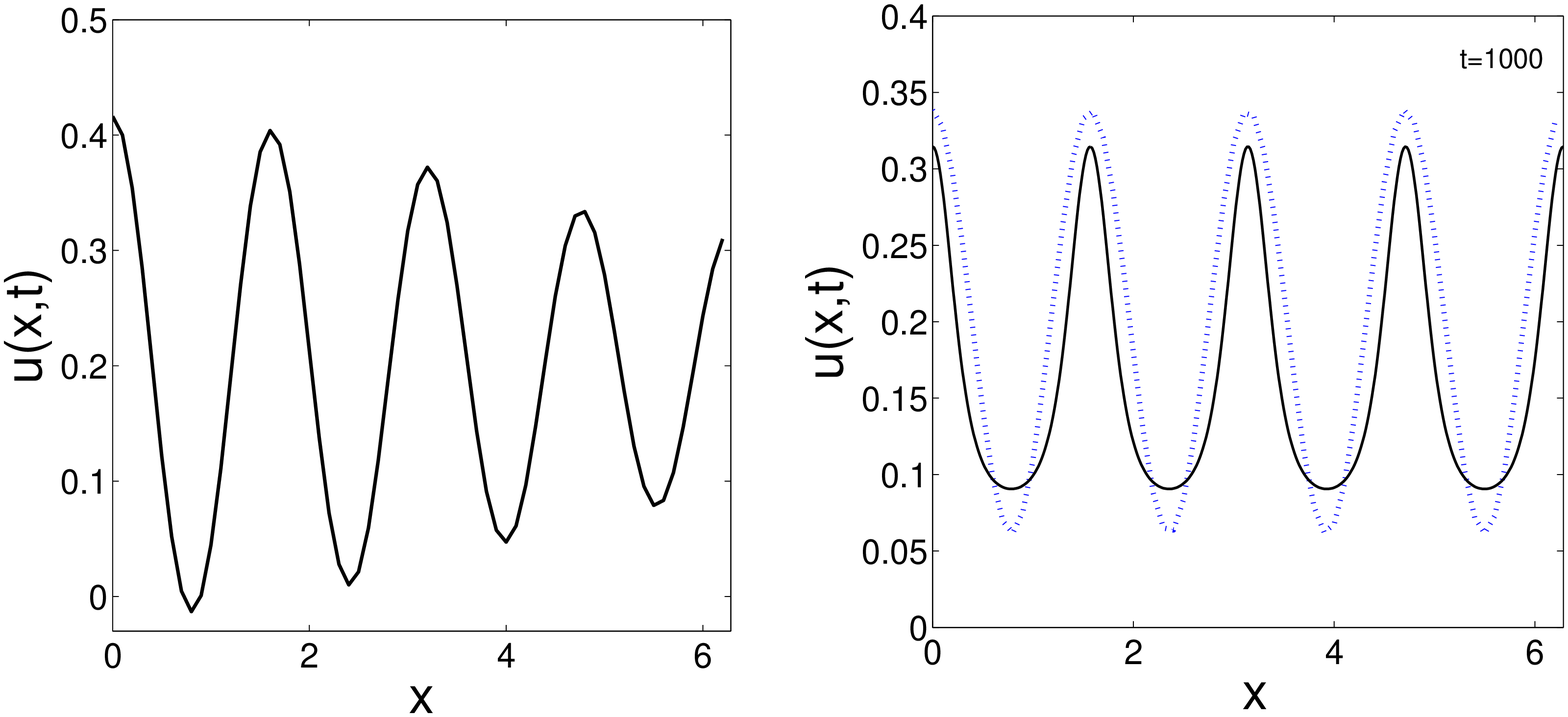}
\caption{Left:Initial condition $u=u_{c}+ \varepsilon\bigl(\bar{A_{1}}M_1\cos(k_{1}x)+ \bar{A_{2}}M_2\cos(k_{2}x)\bigr)$ with $k_{1}=4$ and $k_{2}=3.5$, where $\bar{A_{1}}=0.344, \bar{A_{2}}=0.108$ is in the basin of attraction of the equilibrium $(A_{1},0)$. Right: The comparison between the numerical solution of system \eqref{eq 01}-\eqref{condition Neumann boundary}(solid line) and the expected solution \eqref{eq 50} (dashed line).}\label{Fig 10}
\end{figure*}

It is seen that, after a long time evolution, the mode with the sufficiently dominant initial data can extinguish the other one. Thus, nonlinear effect makes the initial data play a critical role in the competition of unstable modes.

\section{Conclutions}\label{section_5}
In this paper we have analyzed the pattern formation of a volume-filling chemotaxis model with the logistic cell growth. We have first established the sufficient conditions of destabilization for the uniform steady-state by local stability analysis. By deriving the cubic and the quintic Stuart-Landau equations near the instability threshold  $\chi_c$, we have studied the process of pattern formation for both the supercritical case and the subcritical case. For the case of multiple unstable modes, a competition model of two unstable modes has been obtained. Moreover, for the subcritical case, we have verify the existence of patterns with large amplitudes for $\chi<\chi_c$. This gives an affirmative answer to the open question proposed in \cite{MHTT2014}.

All theoretical results have been tested against numerical solutions showing excellent qualitative and good quantitative agreement.

Inspired by the obtained results, we will further study the large amplitude patterns of the system \eqref{eq 01}, such as their structure and stability. In addition, the research idea in this paper can be applied to other biological models, such as the chemotaxis models with different kinetic terms and chemotactic flux introduced by Wang and Xu \cite{WX2013}.

\begin{appendix}
\titleformat{\section}{\Large\bfseries}{Appendix\ \thesection} {1em}{}

\section{Derivation of coefficient equations}\label{appendix coefficient-eqs}

Let
    $$\mathbf{w}=\begin{pmatrix}W^{(1)}\\ W^{(2)}\end{pmatrix} =\begin{pmatrix} u-\bar{u}\\ v-\bar{v}\end{pmatrix},$$
 and then  \eqref{eq 01}  can be rewritten as
\begin{equation}\label{eq A20}
    \mathbf{w_{t}}=\mathcal{L}^{\chi}\mathbf{w} +\mathcal{N}\mathcal{L}^{\chi}\mathbf{w},
\end{equation}
where
    $$\mathcal{L}^{\chi}=K+D^{\chi}\nabla^2,$$
    $$\mathcal{N}\mathcal{L}^{\chi}\mathbf{w}= \begin{pmatrix} -\chi\nabla\bigl((1-2u_{c})W^{(1)} \nabla W^{(2)}-(W^{(1)})^{2}\nabla W^{(2)}\bigr)-\frac\mu{u_{c}} {(W^{(1)})^{2}} \\ 0 \end{pmatrix}.$$
We expand $\chi,\mathbf{w}$ and $t$ as follows
\begin{equation}\label{eq A21}
    \begin{cases}
    \chi=\chi_{c} +\varepsilon\chi_{1} +\varepsilon^{2}\chi_{2} +\varepsilon^{3}\chi_{3} +\varepsilon^{4}\chi_{4} +\varepsilon^{5}\chi_5 +\cdots,\\
    \mathbf{w}=\varepsilon\mathbf{w_{1}} +\varepsilon^{2}\mathbf{w_{2}} +\varepsilon^{3}\mathbf{w_{3}} +\varepsilon^{4}\mathbf{w_{4}} +\varepsilon^{5}\mathbf{W_5}+\cdots,\\
    t=t(T_1,T_2,T_3,\cdots),\quad T_i=\varepsilon^i t, i=1,2,\cdots
    \end{cases}
\end{equation}
and obtain
\begin{equation}\label{eq A31}
    L^{\chi}=L^{\chi_c}-\varepsilon M_{\chi_1}\nabla^2-\varepsilon^2 M_{\chi_2}\nabla^2-\cdots,
\end{equation}
where
    $$M_{\chi_i}=\begin{pmatrix}   0& \chi_i(u_{c}-u_{c}^{2}) \\
    0& 0  \end{pmatrix},\quad i=1,2,\cdots.$$
Substituting  \eqref{eq A21} and \eqref{eq A31}  into \eqref{eq A20} and collecting the terms at each order in $\varepsilon$, we obtain the following sequence of coefficient equations

\noindent $O(\varepsilon)$:
\begin{equation}\label{eq A01}
\mathcal{L}^{\chi_{c}}\mathbf{w_{1}}=0,
\end{equation}
$O(\varepsilon^{2})$:
\begin{equation}\label{eq A02}
\mathcal{L}^{\chi_{c}}\mathbf{w_{2}}=\mathbf{F},
\end{equation}
$O(\varepsilon^{3})$:
\begin{equation}\label{eq A03}
\mathcal{L}^{\chi_{c}}\mathbf{w_{3}}=\mathbf{G},
\end{equation}
$O(\varepsilon^{4})$:
\begin{equation}\label{eq A09}
\mathcal{L}^{\chi_{c}}\mathbf{w_{4}}=\mathbf{H},
\end{equation}
$O(\varepsilon^{5})$:
\begin{equation}\label{eq A10}
\mathcal{L}^{\chi_{c}}\mathbf{w_{5}}=\mathbf{P},
\end{equation}
where:
\begin{align}
    \mathbf{F}=&\frac{\partial\mathbf{w_{1}}}{\partial T_{1}}+M_{\chi_1}\nabla^2\mathbf{w_{1}}\nonumber\\ &+\chi_{c}(1-2u_{c})\nabla\begin{pmatrix} W_{1}^{(1)}\nabla W^{(2)}_{1}\\ 0\end{pmatrix} +\frac{\mu}{u_{c}}\begin{pmatrix} (W^{(1)}_{1})^{2}\\ 0\end{pmatrix} ,\label{eq A32}\\
    \mathbf{G}=&\frac{\partial\mathbf{w_{1}}}{\partial T_{2}}+\frac{\partial\mathbf{w_{2}}}{\partial T_{1}}+M_{\chi_1}\nabla^2\mathbf{w_{2}} +M_{\chi_2}\nabla^2\mathbf{w_{1}} +\frac{2\mu}{u_{c}} \begin{pmatrix} W_{1}^{(1)}W_{2}^{(1)}\\ 0\end{pmatrix}\nonumber\\
    &+\chi_{c}(1-2u_{c})\nabla\begin{pmatrix} W_{2}^{(1)}\nabla W^{(2)}_{1} +W_{1}^{(1)}\nabla W^{(2)}_{2}\\ 0\end{pmatrix}\nonumber\\
    &-\chi_{c}\nabla\begin{pmatrix} (W^{(1)}_{1})^{2}\nabla W^{(2)}_{1}\\ 0\end{pmatrix} +\chi_{1}(1-2u_{c})\nabla\begin{pmatrix} W^{(1)}_{1}\nabla W^{(2)}_{1}\\ 0\end{pmatrix},\label{eq A33}\\
    \mathbf{H}=&\frac{\partial\mathbf{w_{2}}}{\partial T_{2}}+\frac{\partial\mathbf{w_{1}}}{\partial T_{3}}+\frac{\partial\mathbf{w_{3}}}{\partial T_{1}}+M_{\chi_1} \nabla^2\mathbf{w_{3}}+M_{\chi_2} \nabla^2\mathbf{w_{2}}+M_{\chi_3} \nabla^2\mathbf{w_{1}}\nonumber\\
    &+\chi_{c}(1-2u_{c})\nabla\begin{pmatrix} W^{(1)}_{1}\nabla W^{(2)}_{3} +W^{(1)}_{3}\nabla W^{(2)}_{1} +W^{(1)}_{2}\nabla W^{(2)}_{2}\\ 0 \end{pmatrix} \nonumber\\
    &+\chi_{1}(1-2u_{c})\nabla\begin{pmatrix}   W^{(1)}_{1}\nabla W^{(2)}_{2} +W^{(1)}_{2}\nabla W^{(2)}_{1}\\ 0  \end{pmatrix}\nonumber\\
    &+\chi_{2}(1-2u_{c})\nabla\begin{pmatrix} W^{(1)}_{1}\nabla W^{(2)}_{1}\\ 0 \end{pmatrix}\nonumber\\
    &-\chi_{c}\nabla\begin{pmatrix} (W^{(1)}_{1})^{2}\nabla W^{(2)}_{2} +2W^{(1)}_{1}W^{(1)}_{2}\nabla W^{(2)}_{1} \\ 0 \end{pmatrix}\nonumber\\
    &-\chi_{1}\nabla \begin{pmatrix} (W^{(1)}_{1})^{2}\nabla W^{(2)}_{1} \\ 0 \end{pmatrix}+\frac{\mu}{u_{c}} \begin{pmatrix} (W^{(1)}_{2})^{2}+2W^{(1)}_{1}W^{(1)}_{3}\\ 0 \end{pmatrix},\label{eq A34}\\
\intertext{and}
    \mathbf{P}=&\frac{\partial\mathbf{w_{1}}}{\partial T_{4}}+\frac{\partial\mathbf{w_{2}}}{\partial T_{3}}+\frac{\partial\mathbf{w_{3}}}{\partial T_{2}}+\frac{\partial\mathbf{w_{4}}}{\partial T_{1}}\nonumber\\
    &+M_{\chi_1}\nabla^2\mathbf{w_{4}} +M_{\chi_2}\nabla^2\mathbf{w_{3}} +M_{\chi_3}\nabla^2\mathbf{w_{2}} +M_{\chi_4}\nabla^2\mathbf{w_{1}}\nonumber\\
    &+\chi_{c}(1-2u_{c})\nabla\begin{pmatrix} W^{(1)}_{1}\nabla W^{(2)}_{4} +W^{(1)}_{4}\nabla W^{(2)}_{1}\\ 0 \end{pmatrix}\nonumber\\
    &+\chi_{c}(1-2u_{c})\nabla \begin{pmatrix} W^{(1)}_{2}\nabla W^{(2)}_{3} +W^{(1)}_{3}\nabla W^{(2)}_{2} \\ 0 \end{pmatrix}\nonumber\\
    & +\chi_{1}(1-2u_{c})\nabla\begin{pmatrix} W^{(1)}_{2}\nabla W^{(2)}_{2} +W^{(1)}_{1}\nabla W^{(2)}_{3} +W^{(1)}_{3}\nabla W^{(2)}_{1} \\ 0 \end{pmatrix}\nonumber\\
    &+\chi_{2}(1-2u_{c})\nabla\begin{pmatrix}   W^{(1)}_{1}\nabla W^{(2)}_{2}+W^{(1)}_{2}\nabla W^{(2)}_{1}\\ 0 \end{pmatrix}\nonumber\\
    &+\chi_{3}(1-2u_{c})\nabla\begin{pmatrix}   W^{(1)}_{1}\nabla W^{(2)}_{1}\\ 0  \end{pmatrix}\nonumber\\
    &-\chi_{c} \nabla\begin{pmatrix}   2W^{(1)}_{1}W^{(1)}_{3}\nabla W^{(2)}_{1} +2W^{(1)}_{1}W^{(1)}_{2}\nabla W^{(2)}_{2}\\ 0 \end{pmatrix}\nonumber\\
    &-\chi_{c} \nabla\begin{pmatrix}   (W^{(1)}_{1})^{2}\nabla W^{(2)}_{3} +(W^{(1)}_{2})^{2}\nabla W^{(2)}_{1}\\ 0  \end{pmatrix}\nonumber\\
    & -\chi_{1}\nabla\begin{pmatrix} 2W^{(1)}_{1}W^{(1)}_{2}\nabla W^{(2)}_{1} +(W^{(1)}_{1})^{2}\nabla W^{(2)}_{2} \\ 0 \end{pmatrix}\nonumber\\
    &-\chi_{2}\nabla\begin{pmatrix} (W^{(1)}_{1})^{2}\nabla W^{(2)}_{1} \\ 0 \end{pmatrix} +\frac{2\mu}{u_{c}} \begin{pmatrix}   W^{(1)}_{2}W^{(1)}_{3} +W^{(1)}_{1}W^{(1)}_{4}\\ 0 \end{pmatrix}.\label{eq A35}
\end{align}

\section{The cubic Stuart-Landau equation}\label{appendix Stuart-Landau}

Substituting $\mathbf{w_1}$ into  \eqref{eq A32} and combining $\chi_1=\frac{\partial\mathbf{w}}{\partial T_1}=0$ lead to
\begin{equation}\label{eq A04}
\begin{split}
    \mathbf{F}
    \triangleq &\frac 14 A^2\sum_{i=0,2}F_i\cos(ik_cx),
\end{split}\end{equation}
where
    $$F_i=\begin{pmatrix} \frac{2\mu M^2}{u_c}\\ 0\end{pmatrix}-i^2 K_c^2 \begin{pmatrix}M\chi_c(1-2u_c)\\ 0\end{pmatrix},\quad i=0,2.$$
The solution of \eqref{eq A02} with Neumann boundary condition is given by
\begin{equation}\label{eq A05}
    \mathbf{w_{2}}=A^{2}\mathbf{w_{20}} +A^{2}\mathbf{w_{22}}\cos(2k_{c}x).
\end{equation}
Substituting \eqref{eq A05} into \eqref{eq A02}, we have
\begin{equation}\label{eq A06}
    L_i \mathbf{w_{2i}}=\frac 14 F_i,\quad L_i=K-i^2K_c^2 D^{\chi_c},\ i=0,2.
\end{equation}
Further calculations deduce that $\mathbf{w_{20}},\ \mathbf{w_{22}}$ are given as in \eqref{eq 34}.

Substitute $\mathbf{w_1}$, $\mathbf{w_{2}}$ and $\chi_1= \frac{\partial\mathbf{w}}{\partial T_1}=0$ into \eqref{eq A33} and one have
\begin{equation}\label{eq A07}
    \mathbf{G}=\left(\frac{d A}{d T_{2}}\rho -\mathbf{G_{11}}A+\mathbf{G_{13}}A^{3}\right)\cos(k_{c}x) +A^{3}\mathbf{G_3}\cos(3k_{c}x),
\end{equation}
where
\begin{align*}
    \mathbf{G_{11}}=&\begin{pmatrix} 0& k_{c}^{2}\chi_{2}(u_{c}-u_{c}^{2})\\ 0& 0 \end{pmatrix}\rho =\begin{pmatrix} k_{c}^{2}\chi_{2}(u_{c}-u_{c}^{2})\\ 0 \end{pmatrix},\\
    \mathbf{G_{13}}=&\begin{pmatrix} \chi_{c}(2u_{c}-1)k_{c}^{2}\bigl(W^{(1)}_{20}+MW^{(2)}_{22} -\frac 12 W^{(1)}_{22}\bigr)\\ 0  \end{pmatrix}\\
    &+\begin{pmatrix}  \frac{1}{4}\chi_{c}M^2k_{c}^{2} +\frac{\mu M}{u_{c}}(2W^{(1)}_{20}+W^{(1)}_{22})\\ 0 \end{pmatrix},\\
\intertext{and}
    \mathbf{G_3}=&\begin{pmatrix} \chi_{c}(2u_{c}-1)k_{c}^{2}\bigl(3MW^{(2)}_{22} +\frac{3}{2}W^{(1)}_{22}\bigr)\\ 0 \end{pmatrix}\\
    & +\begin{pmatrix} \frac{3}{4}\chi_{c}M^2k_{c}^{2} +\frac{\mu}{u_{c}}MW^{(1)}_{22})\\ 0\end{pmatrix}.
\end{align*}
By the solvability condition  $\langle \mathbf{G},\mathbf{w}^*\rangle=0$ (the form of  $\mathbf{w}^*$ can be seen in\eqref{eq 11}), we obtain the cubic Stuart-Landau  equation \eqref{eq 13} of amplitude $A$. 

\section{The quintic Stuart-Landau equation}\label{appendix Stuart-Landau 5}

Combining with \eqref{eq A07}, the solution of \eqref{eq A03} with Neumann boundary condition is given by
\begin{equation}\label{eq A22}
    \mathbf{w_{3}}=(A\mathbf{w_{31}} +A^{3}\mathbf{w_{32}})\cos(k_{c}x)  +A^{3}\mathbf{w_{33}}\cos(3k_{c}x),
\end{equation}
where $\mathbf{w_{31}},\ \mathbf{w_{32}}$ and $\mathbf{w_{33}}$ satisfy
\begin{equation}\label{eq A23}
    \begin{cases}
        L_1\mathbf{w_{31}}=\sigma\rho-G_{11},\\
        L_1\mathbf{w_{32}}=-L\rho+G_{13},\\
        L_3\mathbf{w_{33}}=G_3
    \end{cases}
\end{equation}
and $L_i,i=1,3$ are given in \eqref{eq A06}.

Substitute  $\mathbf{w_1,\ w_2}$ and $\mathbf{w_3}$  into  \eqref{eq A34}, take $\chi_1= \frac{\partial\mathbf{w}}{\partial T_1}=0$, and then
\begin{align*}
    \mathbf{H}=&\left(\mathbf{w}_{20}\frac{\partial A^2}{\partial T_2}+H_{02}A^2+H_{04}A^4\right) +\left(\frac{\partial A}{\partial T_3}-k_c^2M_{\chi_3}A\right)\rho\cos(k_cx)\\
    &+\left(\mathbf{w}_{22}\frac{\partial A^2}{\partial T_2}+H_{22}A^2+H_{24}A^4\right)\cos(2k_cx)+H_4A^4\cos{(4k_cx)},
\end{align*}
where
    \begin{align*}
        H_{02}=&\frac{\mu M}{u_c} \begin{pmatrix}W_{31}^{(1)}\\ 0 \end{pmatrix},\\ H_{04}=&\frac{\mu}{u_c}\begin{pmatrix} (W_{20}^{(1)})^2 +W_{32}^{(1)}M +\frac 12(W_{22}^{(1)})^2\\ 0 \end{pmatrix},\\
        H_{22}=&-4k_c^2M_{\chi_2}W_{22} +\frac{\mu M}{u_c}\begin{pmatrix} W_{31}^{(1)}\\ 0 \end{pmatrix}\\
        &+(2u_c-1)k_c^2\begin{pmatrix} \chi_c  W_{31}^{(1)} +\chi_c M W_{31}^{(2)}+\chi_2M\\ 0\end{pmatrix},\\
        H_{24}=&\chi_c(2u_c-1)k_c^2\begin{pmatrix} W_{32}^{(2)}M +3W_{33}^{(2)}M +W_{32}^{(1)}\\ 0 \end{pmatrix}\\
        &+\chi_c(2u_c-1)k_c^2\begin{pmatrix} 4W_{20}^{(1)}W_{22}^{(2)}-W_{33}^{(1)}\\ 0 \end{pmatrix}\\
        &+2\chi_c k_c^2\begin{pmatrix} W_{22}^{(2)}M^2+MW_{20}^{(1)}\\ 0 \end{pmatrix}\\
        &+\frac\mu{u_c}\begin{pmatrix} 2W_{20}^{(1)}W_{22}^{(1)} +W_{32}^{(1)}M+W_{33}^{(1)}M\\ 0 \end{pmatrix},\\
    \intertext{and}
        H_{4}=&\chi_c(2u_c-1)k_c^2\begin{pmatrix} 6W_{33}^{(2)}M +2W_{33}^{(1)} +4W_{22}^{(1)}W_{22}^{(2)}\\ 0 \end{pmatrix}\\
        &+2\chi_ck_c^2\begin{pmatrix} W_{22}^{(2)}M^2 +MW_{22}^{(1)}\\ 0 \end{pmatrix}\\
        &+\frac\mu{u_c}\begin{pmatrix}\frac 12(W_{22}^{(1)})^2+W_{33}^{(1)}M\\ 0 \end{pmatrix}.
    \end{align*}
By the solvability condition  $\langle \mathbf{H},\mathbf{w}^*\rangle=0$, we obtain
    $$\frac{\partial A}{\partial T_3}=\gamma A\quad \text{with } \gamma=\frac{k_c^2\chi_3(u_c-u_c^2)M^*}{1+MM^*}>0.$$
And then, $A$ can not cpature the evolution of the amplitude. So, we choose  $\frac{\partial\mathbf{w}}{\partial T_3}=\chi_3=0$ simply to satisfy the solvability condition. Notice of \eqref{eq 13}, the expression $\mathbf{H}$  can be rewritten as
\begin{align*}
    \mathbf{H}=&\left((2\sigma \mathbf{w}_{20}+H_{02})A^2 +(H_{04}-2L\mathbf{w}_{20})A^4\right)\\
    &+\left((2\sigma \mathbf{w}_{22}+H_{22})A^2 +(H_{24}-2L\mathbf{w}_{22})A^4\right)\cos(2k_cx)\\
    &+H_4A^4\cos{(4k_cx)}.
\end{align*}
It is easy to see that the solution of \eqref{eq A09} with Neumann boundary condition is
\begin{equation}\label{eq A24}
    \mathbf{w_{4}}=A^{2}\mathbf{w_{40}} +A^{4}\mathbf{w}_{41}+(A^{2}\mathbf{w_{42}} +A^{4}\mathbf{w}_{43})\cos(2k_{c}x) +A^{4}\mathbf{w}_{44}\cos(4k_{c}x),
\end{equation}
where $\mathbf{w_{40}},\ \mathbf{w}_{41},\ \mathbf{w_{42}},\ \mathbf{w}_{43}$ and $\mathbf{w}_{44}$   satisfy
\begin{equation}\label{eq A25}
    \begin{cases}
        L_0\mathbf{w}_{40}=2\sigma \mathbf{w}_{20}+H_{02},\\
        L_0\mathbf{w}_{41}=-2L\mathbf{w}_{20}+H_{04},\\
        L_2\mathbf{w}_{42}=2\sigma \mathbf{w}_{22}+H_{22},\\
        L_2\mathbf{w}_{43}=-2L\mathbf{w}_{22}+H_{24},\\
        L_4\mathbf{w}_{44}=H_4
    \end{cases}
\end{equation}
and the definition of  $L_i,i=0,2,4$ are given in \eqref{eq A06}.

Substituting $\mathbf{w_{1}}, \mathbf{w_{2}}, \mathbf{w_{3}}, \mathbf{w_{4}}$ into  \eqref{eq A35} and combining $\chi_1=\chi_3=\frac{\partial\mathbf{w}}{\partial T_1}=\frac{\partial\mathbf{w}}{\partial T_3}=0$ lead to
\begin{align*}
    \mathbf{P}=&\biggl(\frac{\partial A}{\partial T_4}\rho+\frac{\partial A}{\partial T_2}\mathbf{w_{31}}+3A^{2}\frac{\partial A}{\partial T_2}\mathbf{w_{32}}\\
    &\quad -A\mathbf{P_{11}}+A^{3}\mathbf{P_{13}} +A^{5}\mathbf{P_{15}}\biggr)\cos(k_{c}x)+\mathbf{P^{*}},
\end{align*}
where
\begin{align*}
    \mathbf{P_{11}}=&k_c^2\left(M_{\chi_2} \mathbf{w}_{31}+M_{\chi_4}\rho\right),\\
    \mathbf{P_{13}}=&-k_c^2 M_{\chi_2}\mathbf{w}_{32}\\
    &+(2u_{c}-1)\chi_{c}k_{c}^{2}\begin{pmatrix} W_{40}^{(1)} +W_{20}^{(1)}W_{31}^{(2)} -\frac{1}{2}W_{22}^{(1)}W_{31}^{(2)}\\ 0\end{pmatrix}\\
    &+(2u_{c}-1)\chi_{c}k_{c}^{2}\begin{pmatrix} W_{31}^{(1)}W_{22}^{(2)}+W_{42}^{(2)}M -\frac{1}{2}W_{42}^{(1)}\\ 0\end{pmatrix}\\
    &+(2u_{c}-1)\chi_{2}k_{c}^{2}\begin{pmatrix} W_{22}^{(2)}M +W_{20}^{(1)} -\frac{1}{2}W_{22}^{(1)}\\ 0 \end{pmatrix}\\
    &+\frac{1}{4}\chi_{c}k_{c}^{2}\begin{pmatrix} 2MW_{31}^{(1)} +W_{31}^{(2)}M^2\\ 0 \end{pmatrix} +\frac{1}{4}\chi_{2}k_{c}^{2}M^2 \begin{pmatrix} 1\\ 0 \end{pmatrix}\\
    &+\frac{2\mu}{u_{c}}\begin{pmatrix} W_{40}^{(1)}M +\frac{1}{2}W_{42}^{(1)}M +W_{20}^{(1)}W_{31}^{(1)} +\frac{1}{2}W_{22}^{(1)}W_{31}^{(1)}\\ 0 \end{pmatrix},\\
    \mathbf{P_{15}}=&(2u_{c}-1)\chi_{c}k_{c}^{2} \begin{pmatrix} W_{41}^{(1)} +W_{20}^{(1)}W_{32}^{(2)} -\frac{1}{2}W_{22}^{(1)}W_{32}^{(2)} +W_{32}^{(1)}W_{22}^{(2)}\\ 0 \end{pmatrix}\\
    &+(2u_{c}-1)\chi_{c}k_{c}^{2}\begin{pmatrix} \frac{3}{2}W_{22}^{(1)}W_{33}^{(2)} -W_{33}^{(1)}W_{22}^{(2)} +W_{43}^{(2)}M
    -\frac{1}{2}W_{43}^{(1)}\\ 0\end{pmatrix}\\
    &+\chi_{c}k_{c}^{2}\begin{pmatrix} \frac{1}{2}MW_{32}^{(1)} -\frac{1}{2}MW_{33}^{(1)} +2W_{20}^{(1)}W_{22}^{(2)}M +\frac{1}{4}W_{32}^{(2)}M^2 \\ 0 \end{pmatrix}\\
    &+\chi_{c}k_{c}^{2}\begin{pmatrix} \frac{3}{4}W_{33}^{(2)}M^2 +(W_{20}^{(1)})^{2} +\frac{1}{2}(W_{22}^{(1)})^{2} -W_{20}^{(1)}W_{22}^{(1)} \\ 0 \end{pmatrix}\\
    &+\frac{2\mu}{u_{c}}\begin{pmatrix} W_{41}^{(1)}M +\frac{1}{2}W_{43}^{(1)}M +W_{20}^{(1)}W_{32}^{(1)} +\frac{1}{2}W_{22}^{(1)}W_{32}^{(1)} +\frac{1}{2}W_{22}^{(1)}W_{33}^{(1)}\\ 0 \end{pmatrix}
\end{align*}
and the expression of  $\mathbf{P^*}$  involves all terms orthogonal to  $\mathbf{w}^*$ . By solvability condition  $\langle \mathbf{P}, \mathbf{w}^* \rangle=0$, we obtain
\begin{equation}\label{eq A26}
    \frac{\partial A}{\partial T_{4}} =\widetilde{\sigma}A -\widetilde{L}A^{3} +\widetilde{Q}A^{5},
\end{equation}
where
\begin{equation}\label{eq A27}
    \begin{cases}
    \widetilde{\sigma}=\dfrac{\langle -\sigma\mathbf{w_{31}} +\mathbf{P_{11}}, \mbox{\boldmath$\psi$}\rangle} {\langle\mbox{\boldmath$\rho$}, \mbox{\boldmath$\psi$}\rangle},\\
    \widetilde{L}=\dfrac{\langle3\sigma\mathbf{w_{32}} -L\mathbf{w_{31}}+\mathbf{P_{13}}, \mbox{\boldmath$\psi$}\rangle} {\langle\mbox{\boldmath$\rho$}, \mbox{\boldmath$\psi$}\rangle},\\
    \widetilde{Q}=\dfrac{\langle3L\mathbf{w_{32}} -\mathbf{P_{15}},\mbox{\boldmath$\psi$}\rangle} {\langle\mbox{\boldmath$\rho$}, \mbox{\boldmath$\psi$}\rangle}.
    \end{cases}
\end{equation}
Adding up \eqref{eq A26} to \eqref{eq 13} , we get the quintic Stuart-Landau equation \eqref{eq 18} of amplitude $A$.

\section{Competition model}\label{appendix competing-eqs}

Substituting \eqref{eq 15} into  \eqref{eq A32}, we have
\begin{align*}
    \mathbf{F}=&\sum^{2}_{l=1}\left[\frac{\partial A_{l}}{\partial T_{1}} -k_{l}^{2}\begin{pmatrix} 0 & \chi_{1}(u_{c}-u_{c}^{2})\\ 0 & 0  \end{pmatrix}A_{l}\right]\rho_{l}\cos(k_{l}x)\\
    &+F_{01}A_1^2+F_{02}A_2^2+F_{21}A_1^2\cos(2k_1x) +F_{22}A_2^2\cos(2k_2x)\\
    &+F_{p}A_1A_2\cos(k_1+k_2)x +F_{m}A_1A_2\cos(k_1-k_2)x,
\end{align*}
where
\begin{align*}
    F_{0i}=&\frac\mu{2u_c}\begin{pmatrix} M_i^2 \\ 0 \end{pmatrix}, \quad F_{2i}=\frac\mu{2u_c}\begin{pmatrix} M_i^2 \\ 0\end{pmatrix} -k_i^2\chi_c(1-2u_c)\begin{pmatrix} M_i \\ 0\end{pmatrix},\ i=1,2,\\
    F_{p}=&\frac{\mu}{u_{c}}M_1M_2\begin{pmatrix} 1\\ 0 \end{pmatrix} +\frac{1}{2}(2u_{c}-1)\chi_{c}\begin{pmatrix} k_{2}^{2}M_{1} +k_{1}^{2}M_{2} +k_{1}k_{2}M_{1} +k_{1}k_{2}M_{2} \\ 0\end{pmatrix},\\
    F_{m}=&\frac{\mu}{u_{c}}M_1M_2\begin{pmatrix} 1\\ 0 \end{pmatrix} +\frac{1}{2}(2u_{c}-1)\chi_{c}\begin{pmatrix} k_{2}^{2}M_{1} +k_{1}^{2}M_{2} -k_{1}k_{2}M_{1} -k_{1}k_{2}M_{2}\\ 0 \end{pmatrix}.
\end{align*}
Similarly, we take $\chi_{1}= \frac{\partial\mathbf{w}}{\partial T_1}=0$ simply to satisfy the solvability condition. Then,
\begin{equation}\label{eq A13}
    \begin{split}
    \mathbf{F}=&F_{01}A_1^2+F_{02}A_2^2+F_{21}A_1^2\cos(2k_1x) +F_{22}A_2^2\cos(2k_2x)\\
    &+F_{p}A_1A_2\cos(k_1+k_2)x +F_{m}A_1A_2\cos(k_1-k_2)x
    \end{split}
\end{equation}
and the solution of  \eqref{eq A02} with Neumann boundary condition is given by
\begin{equation}\label{eq A14}
\begin{split}
    \mathbf{w_{2}}=&\sum_{l=1}^{2}A_{l}^{2} \sum_{i=0,2}\mathbf{w}_{2i}^{l}\cos(ik_{l}x)\\
    &+A_{1}A_{2}\left(\mathbf{w}_{2p}\cos(k_{1}x +k_{2}x)+\mathbf{w}_{2m}\cos(k_{1}x-k_{2}x)\right),
\end{split}
\end{equation}
where $\mathbf{w}_{20}^{l},\ \mathbf{w}_{22}^{l},\ \mathbf{w}_{2p}$ and $\mathbf{w}_{2m}$ satisfy
\begin{equation*}\begin{cases}
    K\mathbf{w}^{l}_{20}=F_{0l},\ l=1,2,\\ \left(K-4k_{l}^{2}D^{\chi_c}\right)\mathbf{w}^l_{22} =F_{2l},\ l=1,2,\\
    \left(K-(k_1+k_2)^2D^{\chi_x}\right)\mathbf{w}_{2p}=F_p,\\
    \left(K-(k_1-k_2)^2D^{\chi_x}\right)\mathbf{w}_{2m}=F_m.
\end{cases}\end{equation*}

Substituting  \eqref{eq 15} and \eqref{eq A14}  into \eqref{eq A33}, and combining  $\chi_1= \frac{\partial\mathbf{w}}{\partial T_1}=0$, we have
\begin{equation}\label{eq A15}
    \mathbf{G}=\sum^{2}_{l=1}\bigg(\frac{d A_{l}}{d T_2}\rho_{l}-A_{l}\mathbf{G}_{11}^{l} +A_{l}^{3}\mathbf{G}_{13}^{l} +\frac{A_{1}^{2}A_{2}^{2}}{A_{l}} \mathbf{G}_{12}^{l}\bigg)\cos(k_{l}x) +\mathbf{G^*},
\end{equation}
where
\begin{align*}
    \mathbf{G}_{11}^{l}=&M_{\chi_2}k_{l}^{2}\rho_{l},\\
    \mathbf{G}_{13}^{l}=& (2u_{c}-1)\chi_{c}k_{l}^{2} \begin{pmatrix} W_{22}^{l(2)}M_l -\frac{1}{2}W_{22}^{l(1)} +W_{20}^{l(1)}\\ 0 \end{pmatrix}\\
    &+\frac{1}{4}\chi_{c}k_{l}^{2} \begin{pmatrix} M_{l}^2 \\ 0 \end{pmatrix} +\frac{\mu}{u_{c}}M_l \begin{pmatrix} 2W_{20}^{l(1)}+W_{22}^{l(1)}\\                   0 \end{pmatrix} ,\\
    \mathbf{G}_{12}^{1}=& \bigg(\frac{1}{2}(2u_{c}-1)\chi_{c}\big((k_{1}^{2} +k_{1}k_{2})W^{(2)}_{2p}M_2+(k_{1}^{2} -k_{1}k_{2})W^{(2)}_{2m}M_2\\
    &\quad +k_{1}k_{2}(W^{(1)}_{2m}-W^{(1)}_{2p} ) +2W_{20}^{2(1)}k_{1}^{2}\big)+\frac{1}{2}\chi_{c}k_{1}^{2}M_{2}^2\\
    &\quad  +\frac{\mu}{u_{c}}(2W_{20}^{2(1)}M_1+W^{(1)}_{2p}M_2 +W^{(1)}_{2m}M_2)\bigg) \begin{pmatrix} 1\\ 0 \end{pmatrix},\\
    \mathbf{G}_{12}^{2}=& \bigg(\frac{1}{2}(2u_{c}-1)\chi_{c}\big((k_{2}^{2} +k_{1}k_{2})W^{(2)}_{2p}M_1+(k_{2}^{2} -k_{1}k_{2})W^{(2)}_{2m}M_1\\
    &\quad +k_{1}k_{2}(W^{(1)}_{2m}-W^{(1)}_{2p}) +2W_{20}^{1(1)}k_{2}^{2}\big)+\frac{1}{2}\chi_{c}k_{2}^{2}M_{1}^2\\
    &\quad  +\frac{\mu}{u_{c}}(2W_{20}^{1(1)}M_2+W^{(1)}_{2p}M_1 +W^{(1)}_{2m}M_1)\bigg) \begin{pmatrix} 1\\ 0 \end{pmatrix}.
\end{align*}
By the solvability condition, we can obtain the competition model \eqref{eq 22} of two unstable modes,  
where
\begin{gather*}
    \sigma_{l}=\frac{\langle \mathbf{G}_{11}^{l}, \mbox{\boldmath$\psi$}_{l}\rangle} {\langle\mbox{\boldmath$\rho$}_{l}, \mbox{\boldmath$\psi$}_{l}\rangle},\quad
    L_{l}=\frac{\langle \mathbf{G}_{13}^{l}, \mbox{\boldmath$\psi$}_{l}\rangle} {\langle\mbox{\boldmath$\rho$}_{l}, \mbox{\boldmath$\psi$}_{l}\rangle},\quad
    \Omega_{l}=\frac{\langle \mathbf{G}_{12}^{l}, \mbox{\boldmath$\psi$}_{l}\rangle} {\langle\mbox{\boldmath$\rho$}_{l}, \mbox{\boldmath$\psi$}_{l}\rangle},\\
    \mbox{\boldmath$\psi$}_{l}= \begin{pmatrix}  M^{*}_{l}\\1 \end{pmatrix},\quad M^{*}_{l}=\frac{\alpha}{\mu+k_{l}^{2}d_{1}},\ l=1,2.
\end{gather*}

%

\end{appendix}

\end{document}